\documentclass[11pt]{amsart}
\usepackage[T1]{fontenc}
\usepackage[utf8]{inputenc}
\usepackage[margin=3cm,includeheadfoot]{geometry}
\usepackage{mathtools,amssymb}
\usepackage[numbers,sort&compress,square]{natbib}
\usepackage{xcolor}
\usepackage[shortlabels]{enumitem}
\usepackage{aliascnt}
\usepackage[colorlinks=true,citecolor=blue,linkcolor=blue,urlcolor=blue]{hyperref}
\usepackage[nameinlink,noabbrev]{cleveref}

\newtheorem{theorem}{Theorem}[section]
\newaliascnt{lemma}{theorem}
\newtheorem{lemma}[lemma]{Lemma}
\aliascntresetthe{lemma}
\newaliascnt{proposition}{theorem}
\newtheorem{proposition}[proposition]{Proposition}
\aliascntresetthe{proposition}
\newaliascnt{corollary}{theorem}
\newtheorem{corollary}[corollary]{Corollary}
\aliascntresetthe{corollary}
\theoremstyle{remark}
\newaliascnt{remark}{theorem}
\newtheorem{remark}[remark]{Remark}
\aliascntresetthe{remark}
\crefname{theorem}{Theorem}{Theorems}
\crefname{lemma}{Lemma}{Lemmas}
\crefname{proposition}{Proposition}{Propositions}
\crefname{corollary}{Corollary}{Corollaries}
\crefname{remark}{Remark}{Remarks}
\allowdisplaybreaks[4]
\numberwithin{equation}{section}

\DeclarePairedDelimiterXPP\pk[1]{\mathbb P}{\{} {\}}{}{#1}
\DeclarePairedDelimiterXPP\E[1]{\mathbb E}{\{} {\}}{}{#1}
\DeclareMathOperator{\Var}{Var}
\DeclareMathOperator{\Cov}{Cov}
\newcommand{\R}{\mathbb R}
\newcommand{\dd}{\mathop{}\!\mathrm d}
\newcommand{\one}{\mathbf 1}
\newcommand{\X}{\boldsymbol X}
\newcommand{\B}{\boldsymbol B}
\newcommand{\ba}{\boldsymbol a}
\newcommand{\bb}{\boldsymbol b}
\newcommand{\bc}{\boldsymbol c}
\newcommand{\bx}{\boldsymbol x}
\newcommand{\by}{\boldsymbol y}
\newcommand{\bz}{\boldsymbol z}
\newcommand{\ve}{\varepsilon}

\begin{document}
\title[Simultaneous ruin after a deterministic split]
{Simultaneous ruin probabilities for Brownian risks\break
with a deterministic splitting time}
\author{Chao Lu}
\address{Department of Actuarial Science, University of Lausanne,
UNIL-Dorigny, 1015 Lausanne, Switzerland}
\email{chao.lu@unil.ch}
\date{}

\begin{abstract}
We consider two Brownian risk processes which coincide up to a fixed time
and subsequently evolve with correlated increments. For a finite time horizon
after the split, we derive the exact asymptotics of their simultaneous ruin
probability as the initial reserves tend to infinity proportionally.
The terminal correlation separates the two asymptotic regimes for the
ratio to the joint tail at the terminal time. If both coordinates contribute to the leading
Gaussian tail, the ratio is a bivariate Piterbarg-type constant, for which
we give an integral representation and explicit bounds. In the remaining
regimes the ratio equals two, with a Gaussian factor in the endpoint tail
at the transition. We also obtain asymptotics on intervals of length
of order $u^{-2}$ near the terminal time and an exponential limit law
for the first simultaneous ruin time, conditional on ruin by the horizon.
\end{abstract}
\maketitle

\section{\textbf{Introduction}}\label{sec:introduction}

Simultaneous ruin in multivariate Brownian risk models describes the event
that several reserve processes become negative at the same time.
For two Brownian risks with constant correlation, exact finite-horizon
asymptotics were obtained in \cite{DebHasMic2020}. Related results concern
Parisian ruin \cite{Kri2022}, multivariate failure events
\cite{DebHasKri2022}, and hitting a cone over an infinite time horizon
\cite{DebHasJiRol2018}. Background on the connection between ruin and
queueing models can be found in \cite{DebMan2015}.

Exact asymptotics for vector-valued Gaussian processes with independent
coordinates are established in \cite{DebHasJiTab2015}; dependent
coordinates are treated in \cite{DebHasWang2020}. Extensions to locally
additive Gaussian fields are developed in \cite{IevKri2025}.
Sojourns of stationary vector-valued Gaussian fields and more general
path functionals of vector-valued Gaussian processes are studied in
\cite{DebHasMic2025Sojourn} and \cite{IevNovSha2026}, respectively.
The main theorems in these works concern centred Gaussian processes or
fields without an added deterministic trend. Consequently, they do not
directly cover the moving thresholds $u+c_1t$ and $au+c_2t$ considered
below; their application would require an additional argument
incorporating the trend.

In this paper the two risks share a common Brownian path until a
deterministic splitting time. After the split, their increments have a
constant correlation, whereas the correlation between the risk levels
depends on time. This distinction enters the asymptotic constant through
both the covariance at the terminal time and the covariance of the
increments near that time.

Let $W_0,W_1,W_2$ be independent standard Brownian motions. Throughout the
paper assume that
\begin{equation}\label{eq:assumptions}
 0<\gamma<T,\qquad -1<\rho<1,\qquad c_1,c_2>0,\qquad a\leq1.
\end{equation}
Put $W_1^*(t)=\rho W_1(t)+\sqrt{1-\rho^2}\,W_2(t)$ and define
\begin{equation}\label{eq:model}
 \begin{split}
 X_1(t)&=W_0(t\wedge\gamma)+W_1((t-\gamma)_+),\\
 X_2(t)&=W_0(t\wedge\gamma)+W_1^*((t-\gamma)_+),\qquad t\geq0.
 \end{split}
\end{equation}
Each coordinate is a standard Brownian motion. The vector process
$\X=(X_1,X_2)^\top$ has independent increments, and
\begin{equation}\label{eq:covariance}
 \Cov\bigl(X_1(s),X_2(t)\bigr)
 =(s\wedge t)\wedge\gamma+\rho\bigl((s\wedge t)-\gamma\bigr)_+.
\end{equation}
In particular, for $t>0$ its correlation is
\begin{equation}\label{eq:correlation}
 r(t)=
 \begin{cases}
 1,&0<t\leq\gamma,\\
 \rho+(1-\rho)\gamma/t,&t>\gamma.
 \end{cases}
\end{equation}

For a compact set $I\subset[0,T]$, write
\begin{equation}\label{eq:ruin}
 p_I(u)=\pk*{\exists t\in I\colon
 X_1(t)-c_1t>u,\ X_2(t)-c_2t>au},\qquad u>0,
\end{equation}
and put $q_t(u)=p_{\{t\}}(u)$. We suppress the fixed parameter $a$ in this
notation. When $a>0$, the thresholds $u$ and $au$ represent positive initial
reserves. Allowing $a\leq0$ also describes the transition to an asymptotically
one-dimensional exceedance event.

We show that $p_{[0,T]}(u)\sim Cq_T(u)$ and identify $C$ in all regimes.
The transition occurs at $a=r(T)$. For $a>r(T)$, the constant is expressed
through a Brownian motion with correlation $\rho$ and a drift
determined by the terminal covariance. For $a\leq r(T)$, it equals two.
Conditional on ruin by $T$, we further show that the time remaining
from first simultaneous ruin to $T$, multiplied by $u^2$, converges to
an exponential distribution with an explicit rate.
The proofs use conditioning at the left endpoint of a short terminal
interval, Brownian scaling, and a stopping-time bound. The main results
are given in \cref{sec:results}, and their proofs in \cref{sec:proofs}.

\section{\textbf{Main Results}}\label{sec:results} 

Write $r=r(T)$ and introduce
\begin{equation}\label{eq:matrices}
 R_\rho=\begin{pmatrix}1&\rho\\\rho&1\end{pmatrix},\qquad
 \Sigma(t)=t\begin{pmatrix}1&r(t)\\r(t)&1\end{pmatrix},\quad t>\gamma.
\end{equation}
Thus $\Sigma'(t)=R_\rho$. Let $\ba=(1,a)^\top$, $\bc=(c_1,c_2)^\top$, and set
\begin{equation}\label{eq:b}
 \bb=\Sigma(T)^{-1}\ba,\qquad
 b_1=\frac{1-ar}{T(1-r^2)},\qquad
 b_2=\frac{a-r}{T(1-r^2)}.
\end{equation}
If $a>r$, then $b_1,b_2>0$. Define also
\begin{equation}\label{eq:theta}
 \theta=\bb^\top R_\rho\bb=b_1^2+2\rho b_1b_2+b_2^2>0.
\end{equation}
The strict positivity follows from positive definiteness of $R_\rho$.

Let $\B=(B_1,B_2)^\top$ be a Brownian motion with covariance
$\Cov(\B(s),\B(t))=(s\wedge t)R_\rho$. For $a>r$ and $L\geq0$, put
\begin{equation}\label{eq:constant-finite}
 \mathcal C_L=
 \int_{\R^2}\pk*{\exists t\in[0,L]\colon
 \B(t)-R_\rho\bb\,t>\bx}
 e^{\bb^\top\bx}\dd\bx,
\end{equation}
where vector inequalities are coordinatewise. We use the notation
\begin{equation}\label{eq:constant-infinite}
 \mathcal C=\int_{\R^2}\pk*{\exists t\geq0\colon
 \B(t)-R_\rho\bb\,t>\bx}
 e^{\bb^\top\bx}\dd\bx.
\end{equation}
In particular, $\mathcal C_0=(b_1b_2)^{-1}$. The dependence of
$\mathcal C_L$ and $\mathcal C$ on $a,\rho,\gamma,T$ is suppressed.

Let $\Phi$, $\Psi=1-\Phi$, and $\varphi$ denote the standard normal
distribution function, survival function, and density, respectively.
For $a\leq r$, set
\begin{equation}\label{eq:h}
 h_a=
 \begin{cases}
 1,&a<r,\\[1mm]
 \displaystyle\Phi\!\left(\frac{(rc_1-c_2)\sqrt T}{\sqrt{1-r^2}}\right),&a=r.
 \end{cases}
\end{equation}
Notice that $h_a>0$, including when $r<0$.

\begin{theorem}\label{thm:main}
Under \eqref{eq:assumptions}, as $u\to\infty$,
\begin{equation}\label{eq:main}
 p_{[0,T]}(u)\sim Cq_T(u),\qquad
 C=
 \begin{cases}
 b_1b_2\mathcal C,&a>r,\\
 2,&a\leq r.
 \end{cases}
\end{equation}
For $a>r$, the constant $\mathcal C$ is finite and satisfies
\begin{equation}\label{eq:constant-bounds}
 \frac1{b_1b_2}\leq\mathcal C
 \leq\frac1{b_1b_2}\frac{4\pi}{\pi+2\arcsin\rho}.
\end{equation}
Consequently, $C$ is independent of $c_1,c_2$.
\end{theorem}

For an explicit form of the equivalent, put
\begin{equation}\label{eq:exponent}
 \begin{split}
 G&=\frac1{2\pi T\sqrt{1-r^2}},\qquad
 V=\frac{1+a^2-2ar}{T(1-r^2)},\\
 \Lambda&=\frac{2\bigl((1-ar)c_1+(a-r)c_2\bigr)}{1-r^2},\qquad
 \Omega=\frac{T(c_1^2+c_2^2-2rc_1c_2)}{1-r^2}.
 \end{split}
\end{equation}

\begin{corollary}\label{cor:explicit}
Under \eqref{eq:assumptions}, as $u\to\infty$,
\begin{equation}\label{eq:explicit}
 p_{[0,T]}(u)\sim
 \begin{cases}
 \displaystyle G\mathcal C u^{-2}
 \exp\!\left(-\frac{Vu^2+\Lambda u+\Omega}{2}\right),&a>r,\\[2mm]
 \displaystyle\frac{2\sqrt T}{\sqrt{2\pi}}h_a u^{-1}
 \exp\!\left(-\frac{(u+c_1T)^2}{2T}\right),&a\leq r.
 \end{cases}
\end{equation}
\end{corollary}

Define the first simultaneous ruin time by
\begin{equation}\label{eq:ruin-time}
 \tau_u=\inf\{t\geq0\colon
 X_1(t)-c_1t>u,\ X_2(t)-c_2t>au\},
 \qquad \inf\varnothing=\infty.
\end{equation}

\begin{theorem}[Conditional ruin time]\label{thm:ruin-time}
Under \eqref{eq:assumptions}, for every $x\geq0$,
\begin{equation}\label{eq:ruin-time-limit}
 \lim_{u\to\infty}
 \pk*{u^2(T-\tau_u)>x\mid\tau_u\leq T}
 =e^{-\delta x/2},\qquad
 \delta=
 \begin{cases}
 \theta,&a>r,\\
 T^{-2},&a\leq r.
 \end{cases}
\end{equation}
Thus, conditional on $\tau_u\leq T$, the variable $u^2(T-\tau_u)$
converges in distribution to an exponential random variable with rate
$\delta/2$.
\end{theorem}

The next two propositions describe the terminal time scale and will also
be used in the proofs of \cref{thm:main,thm:ruin-time}.

\begin{proposition}\label{prop:window}
Under \eqref{eq:assumptions}, for fixed $0\leq j<k<\infty$,
\begin{equation}\label{eq:window}
 \lim_{u\to\infty}
 \frac{p_{[T-k/u^2,T-j/u^2]}(u)}{q_T(u)}
 =
 \begin{cases}
 b_1b_2e^{-\theta j/2}\mathcal C_{k-j},&a>r,\\[1mm]
 2e^{-j/(2T^2)}\Phi\bigl(\sqrt{k-j}/T\bigr),&a\leq r.
 \end{cases}
\end{equation}
\end{proposition}

\begin{proposition}\label{prop:localization}
Under \eqref{eq:assumptions}, for every $0<\ve<T-\gamma$,
\begin{equation}\label{eq:far}
 p_{[0,T-\ve]}(u)=o(q_T(u)),\qquad u\to\infty.
\end{equation}
Moreover,
\begin{equation}\label{eq:localization}
 \lim_{L\to\infty}\limsup_{u\to\infty}
 \frac{p_{[0,T-L/u^2]}(u)}{q_T(u)}=0.
\end{equation}
\end{proposition}

\begin{remark}[Constant correlation]\label{rem:constant-correlation}
The proofs also apply when $\gamma=0$. Then $r(t)=\rho$ for $t>0$.
For $T=1$, we have
\[
 b_1=\frac{1-a\rho}{1-\rho^2},\qquad
 b_2=\frac{a-\rho}{1-\rho^2},\qquad
 R_\rho\bb=(1,a)^\top.
\]
Thus \cref{thm:main} gives the constant-correlation result of
\cite{DebHasMic2020}. For $\gamma>0$, the transition point is
$r(T)=\rho+(1-\rho)\gamma/T>\rho$, and the drift in
\eqref{eq:constant-infinite} is $R_\rho\Sigma(T)^{-1}\ba$.
\end{remark}

\section{\textbf{Proofs}}\label{sec:proofs}

For $t>\gamma$, let $f_t$ denote the density of $\X(t)$, and let
$\varphi_t$ denote the density of a centered normal variable with variance
$t$. Write $\bz_u(t)=u\ba+t\bc$. In particular,
\begin{equation}\label{eq:density}
 f_T(\bz_u(T))=G\exp\!\left(-\frac{Vu^2+\Lambda u+\Omega}{2}\right).
\end{equation}

\subsection{Endpoint tails}

We use the Gaussian tail asymptotics of
\cite[Theorem~4.1]{HashorvaHuesler2003}. Their result allows the
correlation matrix to vary with the threshold and gives the following
estimates after standardization.

\begin{lemma}\label{lem:endpoint}
Under \eqref{eq:assumptions}, for every fixed $k\geq0$ and
$t_u=T-k/u^2$, as $u\to\infty$,
\begin{equation}\label{eq:endpoint}
 q_{t_u}(u)\sim
 \begin{cases}
 \displaystyle\frac{f_{t_u}(\bz_u(t_u))}{b_1b_2u^2},&a>r,\\[2mm]
 \displaystyle\frac{Th_a}{u}\varphi_{t_u}(u+c_1t_u),&a\leq r.
 \end{cases}
\end{equation}
In particular, taking $k=0$ gives the equivalent for $q_T(u)$.
Moreover,
\begin{equation}\label{eq:endpoint-shift}
 \lim_{u\to\infty}\frac{q_{t_u}(u)}{q_T(u)}
 =\begin{cases}
 e^{-\theta k/2},&a>r,\\
 e^{-k/(2T^2)},&a\leq r.
 \end{cases}
\end{equation}
\end{lemma}

For $a=r$ and $c_2>rc_1$, first change the sign of the second coordinate.
Apply the cited theorem to the resulting standardized vector and take the
complement within $\{X_1(t_u)>u+c_1t_u\}$. This yields the same factor $h_a$ in
\eqref{eq:endpoint}. In the other cases the theorem applies directly.
Here $r(t_u)=r+O(u^{-2})$.

The endpoint ratio \eqref{eq:endpoint-shift} follows from
\eqref{eq:endpoint} and the elementary density ratios
\begin{align} 
 \frac{f_{t_u}(\bz_u(t_u))}{f_T(\bz_u(T))}
 &\longrightarrow e^{-\theta k/2} \quad\text{as }u \to \infty,\label{eq:density-shift}\\
 \frac{\varphi_{t_u}(u+c_1t_u)}{\varphi_T(u+c_1T)}
 &\longrightarrow e^{-k/(2T^2)} \quad\text{as }u \to \infty.\label{eq:scalar-density-shift}
\end{align}
For the first ratio, use
$\left.\frac{\dd}{\dd t}(\ba^\top\Sigma(t)^{-1}\ba)\right|_{t=T}
=-\theta$.

\subsection{The Brownian constant}

For $a>r$ and $s\geq0$, define
\begin{equation}\label{eq:H}
 H_s=\int_{\R^2}\pk*{\exists t\in[0,s]\colon\B(t)<\bx}
 e^{-\bb^\top\bx}\dd\bx.
\end{equation}

\begin{lemma}\label{lem:constant}
If $a>r$, then
\begin{equation}\label{eq:tilt-identity}
 \mathcal C_s=e^{-\theta s/2}H_s,\qquad s\geq0.
\end{equation}
Moreover, $\mathcal C_s\uparrow\mathcal C$ as $s\to\infty$, and
\eqref{eq:constant-bounds} holds.
\end{lemma}

\begin{proof}
The assertion at $s=0$ follows by integration over the appropriate
quadrant. Fix $s>0$ and let
\[
 q_\rho=\pk*{B_1(1)>0,B_2(1)>0}
 =\frac14+\frac{\arcsin\rho}{2\pi}>0.
\]
The strong Markov property at the first entrance into
$\bx+(0,\infty)^2$ gives
\begin{equation}\label{eq:orthant-bound}
 \pk*{\exists t\in[0,s]\colon\B(t)>\bx}
 \leq q_\rho^{-1}\pk*{\B(s)>\bx}.
\end{equation}
Indeed, the crossing event equals the event that this entrance time is
strictly less than $s$. At the entrance time both coordinates are at
least their respective thresholds. Conditional on this time and position,
the probability that both subsequent increments are positive is
$q_\rho$. This also covers an entrance at time zero.

By symmetry and Tonelli's theorem, \eqref{eq:orthant-bound} implies
\begin{equation}\label{eq:H-bound}
 H_s\leq q_\rho^{-1}\int_{\R^2}\pk*{\B(s)>\bx}
 e^{\bb^\top\bx}\dd\bx
 =\frac{e^{\theta s/2}}{q_\rho b_1b_2}<\infty.
\end{equation}

To prove \eqref{eq:tilt-identity}, use symmetry in \eqref{eq:H} and
make the pathwise substitution $\bx=\B(s)+\by$. This gives
\[
 e^{-\theta s/2}H_s=
 \E*{e^{\bb^\top\B(s)-\theta s/2}
 \int_{\R^2}\one_{\{\exists t\in[0,s]\colon
 \B(t)-\B(s)>\by\}}e^{\bb^\top\by}\dd\by}.
\]
Under the probability measure with density
$\exp(\bb^\top\B(s)-\theta s/2)$, the process
$\B(t)-R_\rho\bb\,t$, $0\leq t\leq s$, is a centered Brownian motion
with correlation matrix $R_\rho$. This is the Gaussian exponential
change of measure, see also \cite{Lej2018}.
Time reversal of the increments then shows that
$\{\B(s-t)-\B(s),\ 0\leq t\leq s\}$ has the law of
$\{\B(t)-R_\rho\bb\,t,\ 0\leq t\leq s\}$ under the original measure.
This proves \eqref{eq:tilt-identity}.

The definition of $\mathcal C_s$ and monotone convergence give
$\mathcal C_s\uparrow\mathcal C$. The event at time zero yields
$\mathcal C_s\geq(b_1b_2)^{-1}$, and \eqref{eq:H-bound} gives the upper
bound in \eqref{eq:constant-bounds}.
\end{proof}

\subsection{\texorpdfstring{Proof of \cref{prop:window}}{Asymptotics on terminal intervals}}

Fix $0\leq j<k$, put $s=k-j$ and $t_u=T-k/u^2$, and take $u$ large
enough that $t_u>\gamma$. For $\bx=(x_1,x_2)^\top$, define
\begin{equation}\label{eq:D}
 D_{u,s}(\bx)=\pk*{\exists t\in[0,s]\colon
 B_i(t)+c_it/u<x_i,\ i=1,2},\qquad
 D_s(\bx)=\pk*{\exists t\in[0,s]\colon\B(t)<\bx}.
\end{equation}
Since $c_i>0$, the events defining $D_{u,s}(\bx)$ increase to the event
defining $D_s(\bx)$ as $u\to\infty$. To see the converse inclusion in
the limit, choose a witnessing time for the strict inequalities and
then take $u$ large enough to preserve both positive margins. Thus
\begin{equation}\label{eq:D-convergence}
 D_{u,s}(\bx)\longrightarrow D_s(\bx)
 \quad\text{as }u \to \infty \quad\text{for every }\bx\in\R^2.
\end{equation}

\begin{proof}[Proof for $a>r$]  
Conditioning on $\X(t_u)=\bz_u(t_u)+\bx/u$, using independent and stationary
increments, Brownian scaling, and symmetry, yields the exact identity
\begin{equation}\label{eq:window-integral}
 p_{[t_u,t_u+s/u^2]}(u)
 =\frac{f_{t_u}(\bz_u(t_u))}{u^2}
 \int_{\R^2}D_{u,s}(\bx)
 e^{-\bb_u^\top\bx-\bx^\top\Sigma(t_u)^{-1}\bx/(2u^2)}\dd\bx,
\end{equation}
where
\begin{equation}\label{eq:bu}
 \bb_u=\Sigma(t_u)^{-1}\left(\ba+\frac{t_u\bc}{u}\right)
 \longrightarrow\bb \quad\text{as }u \to \infty.
\end{equation}

We verify domination on the whole plane. Put
$M_i=\sup_{0\leq t\leq s}(-B_i(t))$. For every $\beta_1,\beta_2>0$,
Tonelli's theorem gives
\begin{equation}\label{eq:majorant-moment}
 \int_{\R^2}D_s(\bx)e^{-\beta_1x_1-\beta_2x_2}\dd\bx
 \leq\frac{\E*{e^{\beta_1M_1+\beta_2M_2}}}{\beta_1\beta_2}<\infty.
\end{equation}
The last inequality follows from Cauchy--Schwarz and the reflection
principle, since each $M_i$ has the law of $\sqrt s\,|Z|$ for a standard
normal variable $Z$. For large $u$, $b_i/2\leq b_{i,u}\leq3b_i/2$.
Also $D_{u,s}\leq D_s$, and the quadratic term in
\eqref{eq:window-integral} is nonnegative. Hence an integrable majorant
is obtained by summing the four functions in
\eqref{eq:majorant-moment} with $\beta_i\in\{b_i/2,3b_i/2\}$.
Dominated convergence and \eqref{eq:D-convergence} show that the integral
in \eqref{eq:window-integral} converges to $H_s$.

Combining \cref{lem:endpoint} and
\eqref{eq:tilt-identity}, we obtain
\[
 \frac{p_{[T-k/u^2,T-j/u^2]}(u)}{q_T(u)}
 \longrightarrow b_1b_2e^{-\theta k/2}H_s
 =b_1b_2e^{-\theta j/2}\mathcal C_s.\qedhere
\]
\end{proof}

\begin{proof}[Proof for $a\leq r$]
With $\sigma(t)=\sqrt{t(1-r(t)^2)}$, condition on
\[
 X_1(t_u)=u+c_1t_u+x/u,\qquad
 X_2(t_u)=r(t_u)(u+c_1t_u+x/u)+\sigma(t_u)v.
\]
The conditional variable $v$ is standard normal and independent of
$X_1(t_u)$. Define
\begin{equation}\label{eq:yu}
 y_u(x,v)=(r(t_u)-a)u^2+
 \bigl(v\sigma(t_u)-(c_2-r(t_u)c_1)t_u\bigr)u+r(t_u)x.
\end{equation}
As in \eqref{eq:window-integral}, we obtain
\begin{equation}\label{eq:weak-window}
 p_{[t_u,t_u+s/u^2]}(u)
 =\frac{\varphi_{t_u}(u+c_1t_u)}u K_u,
\end{equation}
where
\[
 K_u=\int_{\R^2}D_{u,s}(x,y_u(x,v))
 e^{-\alpha_ux-x^2/(2t_uu^2)}\varphi(v)\dd x\dd v,
 \qquad \alpha_u=t_u^{-1}+c_1/u.
\]
Here and below $D_{u,s}(x,y)$ denotes $D_{u,s}((x,y)^\top)$.

If $a<r$, then $y_u(x,v)\to\infty$ for every fixed $(x,v)$.
If $a=r$, put $v_0=(c_2-rc_1)T/\sigma(T)$. Smoothness of $r$ and
$\sigma$ at $T$ gives
\begin{equation}\label{eq:yu-boundary}
 \begin{split}
 r(t_u)-r&=\frac{(1-\rho)\gamma k}{T^2u^2}+O(u^{-4}),\qquad
 \sigma(t_u)=\sigma(T)+O(u^{-2}),\\
 y_u(x,v)&=u\sigma(T)(v-v_0)
 +\frac{(1-\rho)\gamma k}{T^2}+rx+O(u^{-1}),
 \end{split}
\end{equation}
for fixed $(x,v)$. Consequently, $y_u(x,v)\to\infty$ when $v>v_0$
and $y_u(x,v)\to-\infty$ when $v<v_0$.

Brownian paths are bounded on $[0,s]$. If $y_u\to\infty$, a witnessing
time for $B_1(t)<x$ eventually satisfies both inequalities in
\eqref{eq:D}. If $y_u\to-\infty$, the second inequality eventually
fails at every time. Thus, for almost every $(x,v)$,
\begin{equation}\label{eq:weak-D-limit}
 D_{u,s}(x,y_u(x,v))\longrightarrow
 \pk*{\inf_{0\leq t\leq s}B_1(t)<x}
 \left(\one_{\{a<r\}}+\one_{\{a=r\}}\one_{\{v>v_0\}}\right).
\end{equation}
The exceptional line $v=v_0$ has zero Lebesgue measure.

For large $u$, $1/(2T)\leq\alpha_u\leq2/T$. The integrand defining
$K_u$ is bounded by $m_s(x)\varphi(v)$, where
\[
 m_s(x)=e^{-x/(2T)}\one_{\{x\geq0\}}
 +2e^{-2x/T}\Phi(x/\sqrt s)\one_{\{x<0\}}.
\]
Indeed, for $x<0$ we use $c_1>0$ and the reflection principle to bound
$D_{u,s}(x,y)$ by $2\Phi(x/\sqrt s)$. The function $m_s$ is integrable.
Dominated convergence in \eqref{eq:weak-D-limit} now yields
\begin{align*}
 \lim_{u\to\infty}K_u
 &=h_a\int_\R e^{-x/T}
 \pk*{\inf_{0\leq t\leq s}B_1(t)<x}\dd x\\
 &=Th_a\E*{\exp\!\left(-\frac1T\inf_{0\leq t\leq s}B_1(t)\right)}\\
 &=2Th_a e^{s/(2T^2)}\Phi(\sqrt s/T).
\end{align*}
The last identity follows from
$-\inf_{0\leq t\leq s}B_1(t)\overset d=\sqrt s\,|Z|$.
Combining this limit with \eqref{eq:weak-window} and \cref{lem:endpoint} proves
\eqref{eq:window} for $a\leq r$.
\end{proof} 

\subsection{\texorpdfstring{Proof of \cref{prop:localization}}{Localization}}

We first give the stopping-time bound used to control the interval
between a fixed neighbourhood of $T$ and a terminal interval of length
$L/u^2$. Related finite-horizon bounds appear in \cite{DebHasKri2022}.

\begin{lemma}\label{lem:stopping}
Fix $0<\ve<T-\gamma$, put $d=T-\ve$, and define
\[
 q_\ve=\pk*{\B(\ve)>\ve\bc}>0.
\]
For $d<S\leq T$ and $u>0$,
\begin{equation}\label{eq:stopping}
 p_{[d,S]}(u)\leq q_\ve^{-1}q_S(u).
\end{equation}
\end{lemma}

\begin{proof}
Let $\tau=\inf\{t\geq d\colon\X(t)-t\bc>u\ba\}$, with
$\inf\varnothing=\infty$. By continuity, the crossing event on $[d,S]$
is $\{\tau<S\}$. On this event, $\X(\tau)-\tau\bc\geq u\ba$.
Since $d>\gamma$, the strong Markov property applies to the Brownian
increments after $\tau$, with correlation matrix $R_\rho$.
For every $0<h\leq\ve$, Brownian scaling and $c_1,c_2>0$ give
\[
 \pk*{\B(h)>h\bc}
 =\pk*{\B(1)>\sqrt h\,\bc}
 \geq\pk*{\B(1)>\sqrt\ve\,\bc}=q_\ve.
\]
Conditioning on $\tau$ and the position at $\tau$ therefore yields
$q_S(u)\geq q_\ve\pk*{\tau<S}=q_\ve p_{[d,S]}(u)$.
\end{proof}

\begin{proof}[Proof of \cref{prop:localization}]
Fix $0<\ve<T-\gamma$ and put $d=T-\ve$. Suppose first that $a>r$.
The process $Y(t)=\bb^\top\X(t)$ is a continuous centered Gaussian
process with independent increments. Its variance is
\[
 v(t)=\Var(Y(t))=
 \begin{cases}
 (b_1+b_2)^2t,&0\leq t\leq\gamma,\\
 \theta t+2b_1b_2(1-\rho)\gamma,&\gamma<t\leq T.
 \end{cases}
\]
This function is continuous and strictly increasing, and
$v(T)=\bb^\top\Sigma(T)\bb=V$.
Thus $Y$ has the law of a standard Brownian motion evaluated at $v(t)$.
Since $b_i,c_i>0$ and $\bb^\top\ba=V$, the reflection principle gives
\[
 p_{[0,d]}(u)\leq\pk*{\sup_{0\leq t\leq d}Y(t)>Vu}
 =2\Psi\!\left(\frac{Vu}{\sqrt{v(d)}}\right).
\]
Here $v(d)<V$, so \cref{lem:endpoint} implies
$p_{[0,d]}(u)=o(q_T(u))$. Indeed, the quadratic coefficient in the
logarithm of the ratio is
$-V^2/(2v(d))+V/2<0$.

If $a\leq r$, the same conclusion follows from
\[
 p_{[0,d]}(u)\leq\pk*{\sup_{0\leq t\leq d}X_1(t)>u}
 =2\Psi(u/\sqrt d)
\]
and \cref{lem:endpoint}, because $d<T$ and $h_a>0$.
This proves \eqref{eq:far}.

For fixed $L>0$ and large $u$, let $S_u=T-L/u^2>d$.
By \cref{lem:stopping},
\[
 p_{[0,S_u]}(u)\leq p_{[0,d]}(u)+q_\ve^{-1}q_{S_u}(u).
\]
Using \eqref{eq:far} and \eqref{eq:endpoint-shift}, we obtain
\begin{equation}\label{eq:remainder-bound}
 \limsup_{u\to\infty}\frac{p_{[0,T-L/u^2]}(u)}{q_T(u)}
 \leq q_\ve^{-1}e^{-\delta L/2},
\end{equation}
where $\delta$ is given in \eqref{eq:ruin-time-limit}.
Letting $L\to\infty$ proves \eqref{eq:localization}.
\end{proof}

\subsection{\texorpdfstring{Proofs of \cref{thm:main,cor:explicit}}{Proofs of the main results}}

\begin{proof}[Proof of \cref{thm:main}]
For fixed $L>0$ and sufficiently large $u$,
\[
 p_{[T-L/u^2,T]}(u)\leq p_{[0,T]}(u)
 \leq p_{[T-L/u^2,T]}(u)+p_{[0,T-L/u^2]}(u).
\]
Let
\[
 A_L=
 \begin{cases}
 b_1b_2\mathcal C_L,&a>r,\\
 2\Phi(\sqrt L/T),&a\leq r.
 \end{cases}
\]
By \cref{prop:window} and \eqref{eq:remainder-bound},
\begin{equation}\label{eq:squeeze}
 A_L\leq\liminf_{u\to\infty}\frac{p_{[0,T]}(u)}{q_T(u)}
 \leq\limsup_{u\to\infty}\frac{p_{[0,T]}(u)}{q_T(u)}
 \leq A_L+q_\ve^{-1}e^{-\delta L/2}.
\end{equation}
If $a>r$, then $A_L\to b_1b_2\mathcal C$ by \cref{lem:constant}.
If $a\leq r$, then $A_L\to2$. Letting $L\to\infty$ in
\eqref{eq:squeeze} proves \eqref{eq:main}. The bounds on $\mathcal C$
were proved in \cref{lem:constant}, and its definition depends only on
$a,\rho,\gamma,T$.
\end{proof}

\begin{proof}[Proof of \cref{cor:explicit}]
Combine \cref{thm:main,lem:endpoint} with \eqref{eq:density} and
$\varphi_T(u+c_1T)=(2\pi T)^{-1/2}
\exp(-(u+c_1T)^2/(2T))$.
\end{proof}

\subsection{\texorpdfstring{Proof of \cref{thm:ruin-time}}{Conditional ruin time}}

\begin{proof}
For deterministic $s>\gamma$, continuity implies that $\tau_u=s$
forces equality in at least one of the two ruin inequalities.
The continuous marginal distributions give $\pk*{\tau_u=s}=0$,
while continuity also gives $p_{[0,s]}(u)=\pk*{\tau_u<s}$.
Hence, for fixed $x\geq0$ and all sufficiently large $u$,
\begin{equation}\label{eq:ruin-time-ratio}
 \pk*{u^2(T-\tau_u)>x\mid\tau_u\leq T}
 =\frac{p_{[0,T-x/u^2]}(u)}{p_{[0,T]}(u)}.
\end{equation}

Fix $L>x$. The numerator is bounded by
\begin{align*}
 p_{[T-L/u^2,T-x/u^2]}(u)
 &\leq p_{[0,T-x/u^2]}(u)\\
 &\leq p_{[T-L/u^2,T-x/u^2]}(u)+p_{[0,T-L/u^2]}(u).
\end{align*}
By \cref{prop:window},
\[
 \lim_{u\to\infty}
 \frac{p_{[T-L/u^2,T-x/u^2]}(u)}{q_T(u)}
 =
 \begin{cases}
 b_1b_2e^{-\theta x/2}\mathcal C_{L-x},&a>r,\\
 2e^{-x/(2T^2)}\Phi\bigl(\sqrt{L-x}/T\bigr),&a\leq r.
 \end{cases}
\]
As $L\to\infty$, this expression converges to $Ce^{-\delta x/2}$,
where $C$ is the constant in \cref{thm:main}. The remainder vanishes
in the iterated upper limit by \cref{prop:localization}.
Taking lower and upper limits in the preceding bounds therefore yields
\[
 \lim_{u\to\infty}\frac{p_{[0,T-x/u^2]}(u)}{q_T(u)}
 =Ce^{-\delta x/2}.
\]
Finally, divide by $p_{[0,T]}(u)/q_T(u)\to C$ from
\cref{thm:main} and use \eqref{eq:ruin-time-ratio}.
\end{proof}

\renewcommand{\refname}{

\end{document}